\documentclass[11pt,a4paper]{article}
\usepackage{amsmath}
\usepackage{amsthm}
\usepackage{amssymb}
\usepackage{dsfont}
\usepackage[T1]{fontenc}
\usepackage[french,english,russian]{babel}
\usepackage[utf8]{inputenc}
\usepackage{mathrsfs}
\usepackage{lmodern}
\usepackage{graphicx}
\newtheorem{prop}{Proposition}
\newtheorem{theorem}{Theorem}
\newtheorem*{theorems}{Theorem}

\makeatletter
\def\moverlay{\mathpalette\mov@rlay}
\def\mov@rlay#1#2{\leavevmode\vtop{%
   \baselineskip\z@skip \lineskiplimit-\maxdimen
   \ialign{\hfil$\m@th#1##$\hfil\cr#2\crcr}}}
\newcommand{\charfusion}[3][\mathord]{
    #1{\ifx#1\mathop\vphantom{#2}\fi
        \mathpalette\mov@rlay{#2\cr#3}
      }
    \ifx#1\mathop\expandafter\displaylimits\fi}
\makeatother

\title{Limit laws in the lattice problem. \\
IV. The special case of $\mathbb{Z}^{d}$}
\author{Julien Trevisan}
\begin{document}
\selectlanguage{english}
\maketitle
\bigskip
\selectlanguage{french}
\begin{abstract}
Nous étudions l'erreur du nombre de points du réseau $\mathbb{Z}^{d}$ qui tombent dans un hypercube centré autour de $0$ dilaté et translaté et dont les axes sont parallèles aux axes de coordonnées. Nous montrons que si $t$, le facteur de dilatation, est distribué selon la mesure de probabilité $\frac{1}{T} \rho(\frac{t}{T}) dt$ avec $\rho$ une densité de probabilité sur $[0,1]$, l'erreur, normalisée par $t^{d-1}$, converge en loi lorsque $T \rightarrow \infty$ dans le cas où la translation est de la forme $X=(x, \cdots,x)$ et dans le cas où les coordonnées de $X$ sont indépendantes entre elles, indépendantes de $t$ et distribuées selon la loi uniforme sur $[-\frac{1}{2},\frac{1}{2}]$. Dans les deux cas, on calcule par ailleurs la fonction caractéristique de la loi limite.
\end{abstract}
\selectlanguage{english}
\begin{abstract}
We study the error of the number of points of the lattice $\mathbb{Z}^{d}$ that fall into a dilated and translated hypercube centred around $0$ and whose axis are parallel to the axis of coordinates. We show that if $t$, the factor of dilatation, is distributed according to the probability measure $\frac{1}{T} \rho(\frac{t}{T}) dt$ with $\rho$ being a probability density over $[0,1]$ the error, when normalized by $t^{d-1}$, converges in law when $T \rightarrow \infty$ in the case where the translation is of the form $X=(x,\cdots,x)$ and in the case where the coordinates of $X$ are independent between them, independent from $t$ and distributed according to the uniform law over $[-\frac{1}{2},\frac{1}{2}]$. In both cases, we compute the characteristic function of the limit law.
\end{abstract}
\section{Introduction}
 Let $P$ be a measurable subset of $\mathbb{R}^{d}$ of non-zero finite Lebesgue measure. We want to evaluate the following cardinal number when $t \rightarrow \infty$ : $$ N(tP + X, L) = | (tP + X) \cap L|$$ where $X \in \mathbb{R}^{d}$, $L$ is a lattice of $\mathbb{R}^{d}$ and $t P + X$ denotes the set $P$ dilated by a factor $t$ relatively to $0$ and then translated by the vector $X$.  \\
Under mild regularity conditions on the set $P$, one can show that : 
$$ N(tP + X, L) = t^{d}\frac{\text{Vol}(P)}{\text{Covol}(L)} + o(t^{d}) $$
where $o(f(t))$ denotes a quantity such that, when divided by $f(t)$, it goes to $0$ when $t \rightarrow \infty$ and where $\text{Covol}(L)$ is the volume of a fundamental set of the lattice $L$. \\
We are interested in the error term $$\mathcal{R}(tP + X,L) = N(tP + X, L) - t^{d}\frac{\text{Vol}(P)}{\text{Covol}(L)} \textit{.}$$
In the case where $d=2$ and where $P$ is the unit disk $\mathbb{D}^{2}$, Hardy's conjecture in $\cite{hardy1917average}$ stipulates that we should have for all $\epsilon > 0$, $$\mathcal{R}(t \mathbb{D}^{2}, \mathbb{Z}^{2}) = O(t^{\frac{1}{2}+\epsilon}) $$
where $Y = O(X)$ means that there exists $D > 0$ such that $ |Y| \leqslant D |X|$. \\
One of the result in this direction has been established by Iwaniec and Mozzochi in $\cite{iwaniec1988divisor}$. They have proven that for all $\epsilon > 0$, $$\mathcal{R}(t \mathbb{D}^{2}, \mathbb{Z}^{2}) = O(t^{\frac{7}{11}+\epsilon}) \textit{.}$$
This result has been recently improved by Huxley in \cite{huxley2003exponential}. Indeed, he has proven that : $$\mathcal{R}(t \mathbb{D}^{2}, \mathbb{Z}^{2}) = O(t^{K} \log(t)^{\Lambda}) $$
where $K = \frac{131}{208} $ and $\Lambda = \frac{18627}{8320}$.  \\
In dimension 3, Heath-Brown has proven in $\cite{heath2012lattice}$ that : 
$$\mathcal{R}(t \mathbb{D}^{3}, \mathbb{Z}^{3}) = O(t^{\frac{21}{16}+\epsilon}) \textit{.}$$
These last two results are all based on estimating what are called $\textit{exponential sums}$. Furthermore, they only tackle "deterministic" cases. \\ 
Another result about this problem has been established by Bleher, Cheng, Dyson and Lebowitz in $\cite{bleher1993distribution}$. Let $\rho$ be a probability density on $[0,1]$. They took interest in what is happening when the factor of dilatation $t$ is distributed according to the probability measure $\frac{1}{T} \rho(\frac{t}{T}) dt$ and when $T \rightarrow \infty$. Their result states as following : 
\begin{theorems}[$\cite{bleher1993distribution}$]
There exists a probability density $p$ on $\mathbb{R}$ such that for every piecewise continuous and bounded function $g : \mathbb{R} \longrightarrow \mathbb{R}$, 
$$ \lim_{T \rightarrow \infty} \frac{1}{T} \int_{0}^{T} g(\frac{\mathcal{R}(t \mathbb{D}^{2}, \mathbb{Z}^{2})}{\sqrt{t}}) \rho(\frac{t}{T}) dt = \int_{\mathbb{R}} g(x) p(x) dx \textit{.}$$
Furthermore $p$ can be extended as an analytic function over $\mathbb{C}$ and verifies that for every $\epsilon > 0$, $$p(x) = O(e^{-|x|^{4- \epsilon}})$$ when $x \in \mathbb{R}$ and when $|x| \rightarrow \infty$. 
\end{theorems}
We want to follow this approach on another problem. Namely, let us give $a > 0$ and let's define the following set 
\begin{equation}
\label{eq0}
\mathcal{C}(a) = \{ x=(x_{1},\cdots,x_{d}) \in \mathbb{R}^{d} \text{ } | \text{} \forall i \in [1,d] \textit{, } |x_{i}| \leqslant a \} \textit{.}
\end{equation}
In that case, with $\rho$ being a probability density over $[0,1]$, we want to study the possible convergence in distribution of the quantity $\frac{\mathcal{R}(t \mathcal{C}(a) + X, \mathbb{Z}^{d}}{t^{d-1}}$. We already proved such a result when the dimension $d$ was equal to $2$. Here, we are going to prove the two following theorems (that constitute a generalization of the previous result) : 
\begin{theorem}
\label{thm1}
For all $x \in \mathbb{R}$, when $t \in [0,T]$ is distributed according to the probability density $\frac{1}{T} \rho( \frac{\cdot}{T})$ on $[0,T]$ then, when $T \rightarrow \infty$, $\frac{\mathcal{R}(t \mathcal{C}(a) + X ,\mathbb{Z}^{d})}{t^{d-1}}$ converges in law with $X= (x,\cdots,x)$. Furthermore, the limit law has the following characteristic function 
$$ \varphi(u) = \frac{ \sin(d 2^{d-1}u y) + \sin(d 2^{d-1}u(1-y))}{d 2^{d-1}u } $$
with $y = |t_{2,0} - t_{1,0}|$ where $t_{2,0}$ is the first $t \geqslant 0$ such that $ -t + x \in \mathbb{Z}$ and $t_{1,0}$ is the first $t \geqslant 0$ such that $ t + x \in \mathbb{Z}$. In fact, $y = |1-2\{x \}| $ où $\{x \}$ désigne la partie fractionnaire de $x$. 
\end{theorem}
\begin{theorem}
\label{thm2}
Let's assume that $x_{1},\cdots,x_{d}$ are independent random variables distributed according to the uniform distribution over $[- \frac{1}{2}, \frac{1}{2} ]$. Let's assume also that $t \in [0,T]$ is distributed according to the probability density $\frac{1}{T} \rho( \frac{\cdot}{T})$ on $[0,T]$. Let's suppose that $t$ and $x_{1},\cdots,x_{d}$ are independent betweem them then, when $T \rightarrow \infty$, $\frac{\mathcal{R}(t \mathcal{C}(1) + X ,\mathbb{Z}^{d})}{t^{d-1}}$ converges in distribution with $X= (x_{1},\cdots,x_{d})$. \\
Furthermore, the limit law has the following characteristic function : 
$$\varphi(u) = (2\frac{1-\cos(2^{d-1}u)}{(2^{d-1}u)^{2}})^{d} \textit{.}$$
\end{theorem}
In particular, we see that the normalization in these two cases of the error $\mathcal{R}$ is of order $t^{d-1}$. Furthermore, the two cases studied here are two extreme cases : the case of Theorem $\ref{thm1}$ is a case where all the $x_{i}$ are linked (in fact, they are all equal) whereas the case of Theorem $\ref{thm2}$ is a case where all the $x_{i}$ are independent between them. \\
Before beginning, let's observe that it is enough to prove Theorem $\ref{thm1}$ and Theorem $\ref{thm2}$ in the case where $\rho = \mathbf{1}_{[0,1]}$ (see, for example, the proof of Theorem 4.2 in $\cite{bleher1992distribution}$). So, in the rest of the article, we are going to suppose that $\rho = \mathbf{1}_{[0,1]}$. \\
Furthermore, regarding Theorem $\ref{thm1}$, we can suppose also that $a=1$. Indeed, instead of considering $t$, we can consider the new variable $ \tilde{t} = a t$ that will be distributed according to the probability measure $\frac{1}{a T } \rho(\frac{t}{a T}) dt$. 
In the next section we are going to give a bit of heuristic about Theorem $\ref{thm1}$ and Theorem $\ref{thm2}$.
\section{A bit of heuristic and plan of the paper}
First, let's say that the normalization of $\mathcal{R}(t \mathcal{C}(1) + X ,\mathbb{Z}^{d})$ by $t^{d-1}$ is quite natural. Indeed, to within a multiplicative factor, it corresponds to the surface measure of $\partial(t \mathcal{C}(1) + X)$. \\
This normalization appears when looking at the following expression of $\mathcal{R}(t \mathcal{C}(1) + X ,\mathbb{Z}^{d})$ : 
$$ \frac{\mathcal{R}(t \mathcal{C}(1) + X ,\mathbb{Z}^{d})}{t^{d-1}} = \frac{\sum_{i=1}^{d} (2t)^{i-1} (\lfloor t + x_{i} \rfloor - \lceil -t + x_{i} \rceil +1 - 2t ) \prod_{j=i+1}^{d} (\lfloor t + x_{j} \rfloor - \lceil -t + x_{j} \rceil +1  ) }{t^{d-1}}  $$ 
with $X= (x_{1},\cdots,x_{d})$ (see Proposition $\ref{prop1}$ and Equation $\ref{eq4}$). \\
\\
$\textbf{Plan of the paper.}$ After having proved this expression of $\mathcal{R}(t \mathcal{C}(1) + X ,\mathbb{Z}^{d})$, we show in section 3 that the study of $\frac{\mathcal{R}(t \mathcal{C}(1) + X ,\mathbb{Z}^{d})}{t^{d-1}}$ can be reduced to the study of a simpler quantity which is $\Delta(t,X)$ (see Proposition $\ref{prop2}$ for the definition of $\Delta(t,X))$. \\
Then, in section 4, we give the proof of Theorem $\ref{thm1}$. In fact, the case of Theorem $\ref{thm1}$ corresponds to a case where the expression $\Delta(t,X)$ (and the expression of $\mathcal{R}(t \mathcal{C}(1) + X ,\mathbb{Z}^{d})$) is simpler. This simple expression is used to compute the characteristic function of $\Delta(t,X)$ (see Proposition $\ref{prop4}$). We conclude by using Levy's continuity theorem (see Theorem $\ref{thm3}$). \\
In section 5, we also compute the characteristic function of $\Delta(t,X)$ in the case of Theorem $\ref{thm2}$ (see Proposition $\ref{prop6}$). In fact, the key of Theorem $\ref{thm2}$ is the independence of the variables $t, x_{1}, \cdots, x_{d}$ which enables us to make this computation. Other computations, with other distributions for the variables $x_{1},\cdots,x_{d}$, but always with the independence theorem, could be made. We, again, conclude by using Levy's continuity theorem. \\
\\
The next section is dedicated to reduce $\mathcal{R}(t \mathcal{C}(1) + X ,\mathbb{Z}^{d})$ when we study its asymptotical behaviour. 
\section{Simplification of the study of $\frac{\mathcal{R}(t \mathcal{C}(1) + X ,\mathbb{Z}^{d})}{t^{d-1}}$}
The main object of this section is to prove the following proposition : 
\begin{prop}
\label{prop2}
One has that : 
\begin{equation}
\label{eq5}
 \frac{\mathcal{R}(t \mathcal{C}(1) + X ,\mathbb{Z}^{d})}{t^{d-1}}  -  \Delta(t,X)  \underset{ t \rightarrow \infty}{\rightarrow} 0
\end{equation}
where $\Delta(t,X)$ is defined by 
\begin{equation}
\label{eq6}
\Delta(t,X) = 2^{d-1} \sum_{i=1}^{d}  (\lfloor t + x_{i} \rfloor - \lceil -t + x_{i} \rceil +1 - 2t )  
\end{equation}
with $X=(x_{1},\cdots,x_{d})$ and the convergence in Equation $(\ref{eq5})$ is uniform in $X \in \mathbb{R}^{d}$.
\end{prop}
It is a proposition that enables to do some reduction about the asymptotical study of $\frac{\mathcal{R}(t \mathcal{C}(1) + X ,\mathbb{Z}^{d})}{t^{d-1}}$. The main idea is that, in this case, everything can be computed quite easily, it is only a matter of definitions. 
\subsection{An expression of $\mathcal{R}(t \mathcal{C}(1) + X ,\mathbb{Z}^{d})$}
The main object of this subsection is to prove the following proposition : 
\begin{prop}
\label{prop1}
We have for every $X \in \mathbb{R}^{d}$, for every $t > 0$, that 
\begin{equation}
\label{eq1}
\mathcal{R}(t \mathcal{C}(1) + X ,\mathbb{Z}^{d}) = \prod_{i=1}^{d}(\lfloor t + x_{i} \rfloor - \lceil -t + x_{i} \rceil +1) -  (2t)^{d} 
\end{equation}
where $X=(x_{1},\cdots,x_{d})$.
\end{prop}
The proof is quite straightforward.
\begin{proof}
Let $X=(x_{1},\cdots,x_{d}) \in \mathbb{R}^{d}$. Let $t > 0$. \\
One has that : 
\begin{align}
\label{eq2}
N(t \mathcal{C}(1) + X ,\mathbb{Z}^{d})& = \left( \sum_{ \substack{ (n_{1},\cdots,n_{d}) \in \mathbb{Z}^{d} \\ \forall i \in [1,d] \textit{, } -t+x_{i} \leqslant n_{i}  \leqslant t+x_{i}  }} 1 \right) \nonumber \\
& = \prod_{i=1}^{d}(\lfloor t + x_{i} \rfloor - \lceil -t + x_{i} \rceil +1) \textit{.}
\end{align}
according to Equation $(\ref{eq0})$. \\
Furthermore, one has that : 
\begin{equation}
\label{eq3}
\text{Vol}(t \mathcal{C}(1) + X)= (2 t)^{d}
\end{equation}
So, Equation $(\ref{eq2})$ and Equation $(\ref{eq3})$ and the definition of $\mathcal{R}(t \mathcal{C}(1) + X ,\mathbb{Z}^{d})$ give us Equation $(\ref{eq1})$.
\end{proof}
With Equation $(\ref{eq1})$, one has that : 
\begin{equation}
\label{eq4}
\frac{\mathcal{R}(t \mathcal{C}(1) + X ,\mathbb{Z}^{d})}{t^{d-1}} = \frac{\sum_{i=1}^{d} (2t)^{i-1} (\lfloor t + x_{i} \rfloor - \lceil -t + x_{i} \rceil +1 - 2t ) \prod_{j=i+1}^{d} (\lfloor t + x_{j} \rfloor - \lceil -t + x_{j} \rceil +1  ) }{t^{d-1}}  \textit{.}
\end{equation}
Thanks to this last remark, the asymptotical study of $\frac{\mathcal{R}(t \mathcal{C}(1) + X ,\mathbb{Z}^{d})}{t^{d-1}}$ with $t$ distributed on $[0,T]$ according to the probability measure $\frac{1}{T} \rho(\frac{t}{T}) dt$ is going to be reduced to the study of a simpler quantity. It is the object of the next subsection.
\subsection{Reduction of the study of $\frac{\mathcal{R}(t \mathcal{C}(1) + X ,\mathbb{Z}^{d})}{t^{d-1}}$}
The main object of this subsection is to prove Proposition $\ref{prop2}$.  
The proof is quite straightforward and lie on the definitions of $\lfloor \cdot \rfloor$ and of $\lceil \cdot \rceil$ and on Proposition $\ref{prop1}$.
\begin{proof}[Proof of Proposition $\ref{prop2}$]
For every $t > 0$, for every $x \in \mathbb{R}$, 
\begin{equation}
\label{eq7}
t+x - 1< \lfloor t + x \rfloor \leqslant t+x
\end{equation}
and 
\begin{equation}
\label{eq8}
-t+x \leqslant \lceil -t + x \rceil < -t+x +1 \textit{.}
\end{equation}
From Equation $(\ref{eq7})$ and Equation $(\ref{eq8})$, one gets that : 
\begin{equation}
\label{eq9}
2 t -1 < (\lfloor t + x \rfloor - \lceil -t + x \rceil  +1) \leqslant 2t +1
\end{equation}
So, from this last equation and from Equation $(\ref{eq4})$, one has that, when $t \rightarrow \infty$, 
\begin{equation}
\label{eq10}
|\frac{\mathcal{R}(t \mathcal{C}(1) + X ,\mathbb{Z}^{d})}{t^{d-1}}  -  \Delta(t,X)| \leqslant \sum_{i=1}^{d-1} 2^{i-1} \frac{1}{t^{d-i}} = O(\frac{1}{t}) \textit{.}
\end{equation}
So, one gets the wanted result.
\end{proof}
Thanks to Proposition $\ref{prop2}$, we see that the asymptotical study of $\frac{\mathcal{R}(t \mathcal{C}(1) + X ,\mathbb{Z}^{d})}{t^{d-1}}$ can be reduced to the study of $\Delta(t,X)$. We are going to use this fact in the next two sections.
\section{Proof of Theorem $\ref{thm1}$}
The main object of this section is to prove Theorem $\ref{thm1}$. We are going to use the reduction that was mentioned before (see Proposition $\ref{prop2}$). In the case of Theorem $\ref{thm1}$, the expression of $\Delta(t,X)$ is simple and the proof of Theorem $\ref{thm1}$ is only a matter of computation of a characteristic function. 
\subsection{Reduction of the study of $\frac{\mathcal{R}(t \mathcal{C}(1) + X ,\mathbb{Z}^{d})}{t^{d-1}}$}
The main object of this subsection is to prove the following proposition : 
\begin{prop}
\label{prop3}
For every $x \in \mathbb{R}$, for every $g \in C_{c}(\mathbb{R})$, 
\begin{equation}
\label{eq11}
\int_{t=0}^{T} (g\left( \frac{\mathcal{R}(t \mathcal{C}(1) + X ,\mathbb{Z}^{d})}{t^{d-1}} \right) - g\left( \Delta(t,X) \right)) \frac{1}{T} dt \underset{T \rightarrow \infty}{\rightarrow} 0
\end{equation}
where $X = (x, \cdots,x)$ and where $\Delta(t,X)$ was defined in Proposition $\ref{prop2}$.
\end{prop}
It should be noted in this case that 
\begin{equation}
\label{eq12}
\Delta(t,X) = d 2^{d-1}(\lfloor t + x \rfloor - \lceil -t + x \rceil +1 - 2t) \textit{.} 
\end{equation}
The proof of Proposition $\ref{prop3}$ is quite straightforward and based on Proposition $\ref{prop2}$.
\begin{proof}
One has for every $0 < \kappa < \frac{1}{2}$ : 
\begin{align}
\label{eq13}
&| \int_{t=0}^{T} (g\left( \frac{\mathcal{R}(t \mathcal{C}(1) + X ,\mathbb{Z}^{d})}{t^{d-1}} \right) - g\left( \Delta(t,X) \right)) \frac{1}{T}  dt |  \leqslant   2 \lVert  g \rVert_{\infty} \int_{0}^{\kappa}  dt  \nonumber \\
& + \int_{\kappa T}^{T} |g\left( \frac{\mathcal{R}(t \mathcal{C}(1) + X ,\mathbb{Z}^{d})}{t^{d-1}} \right) - g\left( \Delta(t,X) \right)| \frac{1}{T}  dt
\end{align}
and, because $g \in C_{c}(\mathbb{R})$, it is a uniformly continuous function and so one has, because of Proposition $\ref{prop2}$, that 
\begin{equation}
\label{eq14}
\limsup_{T \rightarrow \infty} \int_{\kappa T}^{T} |g\left( \frac{\mathcal{R}(t \mathcal{C}(1) + X ,\mathbb{Z}^{d})}{t^{d-1}} \right) - g\left( \Delta(t,X) \right)| \frac{1}{T}  dt = 0
\end{equation}
So, Equation $(\ref{eq13})$ and Equation $(\ref{eq14})$ give us that for every $0 < \kappa \leqslant \frac{1}{2}$ : 
\begin{equation}
\label{eq15}
\limsup_{T \rightarrow \infty} | \int_{t=0}^{T} (g\left( \frac{\mathcal{R}(t \mathcal{C}(1) + X ,\mathbb{Z}^{d})}{t^{d-1}} \right) - g\left( \Delta(t,X) \right)) \frac{1}{T}  dt |  \leqslant   2 \lVert  g \rVert_{\infty} \int_{0}^{\kappa}  dt   \textit{.}
\end{equation}
By making $\kappa$ go to $0$, one gets the wanted result.
\end{proof}
In the next subsection, we are going to compute the characteristic function of $\Delta(t,X)$ to within a multiplicative factor.
\subsection{Computation of characteristic function}
Before stating the main proposition of this section, we need to make some observations and put in place some notations.\\
Let's call $t_{1,0} < \cdots < t_{1,l}$ the different times $t \in [0,T]$ such that $t + x \in \mathbb{Z}$. \\
In the same way, let's call $t_{2,0} < \cdots < t_{2,h}$ the different times $t \in [0,T]$ such that $-t + x \in \mathbb{Z}$. \\
Let's observe that for every $i \in \{0, \cdots, l-1 \}$, $t_{1,i+1} - t_{1,i} = 1$ and that for every $j \in \{0, \cdots, h-1 \}$, $t_{2,j+1} - t_{2,j} = 1$. \\
As a consequence, one has necessarily that $t_{2,0} \in [t_{1,0},t_{1,1}[$ or $t_{1,0} \in [t_{2,0},t_{2,1}[$ and $h=l$ or $h=l-1$ or $h=l+1$. \\
Let's set : 
\begin{equation}
\label{eq16}
y = |t_{1,0} - t_{2,0}| 
\end{equation}
and 
\begin{equation}
\label{eq17}
\tilde{\Delta}(t,x) = (\lfloor t + x \rfloor - \lceil -t + x \rceil +1 - 2t) \textit{.}
\end{equation}
By the way, let's remark that for all $x \in \mathbb{R}$ : 
\begin{equation}
\label{eq24}
y = |1 - 2 \{x \}|
\end{equation}
where $\{ x \}$ stands for the fractional part of the real $x$. \\
Then, with these notations, one has that : 
\begin{prop}
\label{prop4}
For every $x \in \mathbb{R}$, one has that the characteristic function $\varphi_{\tilde{\Delta}(\cdot,x)}$ of $\tilde{\Delta}(t,x)$, with $t$ being distributed according to $\frac{1}{T} \mathbf{1}_{[0,T]} dt$ verifies that for every $u \in \mathbb{R}$, 
$$\varphi_{\tilde{\Delta}(\cdot,x)}(u) =  \frac{h}{u T}( \sin(u y) + \sin(u(1-y))) + O(\frac{1}{T}) $$
where the $O$ is uniform in $x \in \mathbb{R}$. \\
As a consequence, when $T \rightarrow \infty$, for every $u \in \mathbb{R}$, one has that : 
$$\varphi_{\tilde{\Delta}(\cdot,x)}(u)  \underset{T \rightarrow \infty}{\rightarrow} \frac{ \sin(u y) + \sin(u(1-y))}{u } \textit{.}$$

\end{prop}
The proof consists basically in cutting the interval $[0,T]$ into subintervals where all the quantities that intervene in the computation can be expressed simply. 
\begin{proof}
By symmetry, we can, and we will, suppose that $t_{2,0} \in [t_{1,0},t_{1,1}[$. Let $u \in \mathbb{R}$. \\
One has then that :
\begin{equation}
\label{eq18}
\mathbb{E}(e^{i u \tilde{\Delta}(t,x)}) = \sum_{i=0}^{h-1} \int_{t_{1,i}}^{t_{2,i}} e^{i u \tilde{\Delta}(t,x)} \frac{1}{T} dt + \int_{t_{2,i}}^{t_{1,i+1}} e^{i u \tilde{\Delta}(t,x)} \frac{1}{T} dt + O(\frac{1}{T})
\end{equation}
where the $O$ corresponds to the rest of the integral that is calculated on a union of two intervals of respective lengths at most $2$. \\
Let $i \in \{0, \cdots, h-1\}$. \\
Then one has : 
\begin{equation}
\label{eq19}
\int_{t_{1,i}}^{t_{2,i}} e^{i u \tilde{\Delta}(t,x)} \frac{1}{T} dt = \int_{t_{1,i}}^{t_{2,i}}e^{i u (t_{1,i} + t_{2,i} - 2t)}  \frac{1}{T} dt
\end{equation}
according to the Equation $(\ref{eq17})$. \\
So, one gets that : 
\begin{equation}
\label{eq20}
\int_{t_{1,i}}^{t_{2,i}} e^{i u \tilde{\Delta}(t,x)} \frac{1}{T} dt = \frac{\sin(u y)}{u T}
\end{equation}
where one conveys that $ \frac{sin(0)}{0} = 1$ and one has this last equation because for all $j \in \{0,\cdots,h-1 \}$, $ y= t_{2,0}-t_{1,0} = t_{2,j} - t_{1,j} $. \\
In a similar way, one gets that : 
\begin{equation}
\label{eq21}
\int_{t_{2,i}}^{t_{1,i+1}} e^{i u \tilde{\Delta}(t,x)} \frac{1}{T} dt =  \frac{\sin(u(1- y))}{u T}
\end{equation}
because $1-y =  t_{1,i+1} - t_{2,i}$. \\
So, with Equation $(\ref{eq18})$, Equation $(\ref{eq20})$ and Equation $(\ref{eq21})$, one gets that : 
\begin{equation}
\label{eq22}
\mathbb{E}(e^{i u \tilde{\Delta}(t,x)}) = \sum_{i=0}^{h-1} \frac{\sin(u y)}{u T} + \frac{\sin(u(1- y))}{u T} + O(\frac{1}{T}) = \frac{h}{u T}( \sin(u y) + \sin(u(1-y))) + O(\frac{1}{T})
\end{equation}
By using the fact that $ \lim_{T \rightarrow \infty} \frac{h}{T} = 1$, one gets from equation $(\ref{eq22})$ that : 
\begin{equation}
\label{eq23}
\mathbb{E}(e^{i u \tilde{\Delta}(t,x)})  \underset{T \rightarrow \infty}{\rightarrow} \frac{ \sin(u y) + \sin(u(1-y))}{u } \textit{.}
\end{equation}
\end{proof}
\subsection{Conclusion}
To conclude the proof of Theorem $\ref{thm2}$, we need to recall the Lévy's continuity theorem.
\begin{theorem}
\label{thm3}
Let's give us $(X_{n})_{n \geqslant 1}$ a sequence of real random variables and let's call $(\phi_{n})_{n \geqslant 1}$ the associated sequence of their characteristic functions.\\
Let's suppose that the sequence $(\varphi_{n})_{n \geqslant 1}$ converges point wisely to some function $\varphi$. \\
Then, it is equivalent to say that there exists $X$ a real random variable such that $(X_{n})$ converges in law towards $X$ and to say that the function $\varphi$ is continuous at the point $t=0$. \\  
Furthermore, if the last condition is realized, $\varphi$ is the characteristic function of such a $X$. 
\end{theorem}
We can now conclude the proof of Theorem $\ref{thm2}$.
\begin{proof}[Proof of Theorem $\ref{thm2}$]
Because of Proposition $\ref{prop3}$, it is enough to study the asymptotic convergence in law, when $T \rightarrow \infty$, of the quantity $\Delta(t,X)$ with $X= (x,\cdots,x)$ and $t$ being distributed according to the density $\frac{1}{T} \mathbf{1}_{[0,T]}(t) dt$. \\
The fact that 
\begin{equation}
\label{eq25}
\Delta(t,X) = d 2^{d-1}
\end{equation}
and Proposition $\ref{prop4}$ and Theorem $\ref{thm3}$ give us that $\Delta(t,X)$ converges in law, when $T \rightarrow \infty$, and the characteristic function of the limit law is given by 
\begin{equation}
\label{eq26}
\varphi(u) = \frac{ \sin(d 2^{d-1}u y) + \sin(d 2^{d-1}u(1-y))}{d 2^{d-1}u }
\end{equation}
with $y = 1 - 2 \{x \}$ according to Equation $(\ref{eq24})$. 
\end{proof}
\section{Proof of Theorem $\ref{thm2}$}
The main object of this section is to prove Theorem $\ref{thm2}$. We are going to use the reduction that was mentioned before (see Proposition $\ref{prop2}$). In the case of Theorem $\ref{thm2}$, the proof is only a matter of computation of the characteristic function of $\Delta(t,X)$ and here it can be easily dealt with thanks to the independence between the $x_{i}$ and thanks to the independence between the $x_{i}$ and $t$. 
\subsection{Reduction of the study of $\frac{\mathcal{R}(t \mathcal{C}(1) + X ,\mathbb{Z}^{d})}{t^{d-1}}$}
The main object of this subsection is the following proposition : 
\begin{prop}
\label{prop5}
One has that : 
\begin{equation}
\label{eq27}
\frac{ \mathcal{R}(t \mathcal{C}(1) + X ,\mathbb{Z}^{d})}{t^{d-1}}  -  \Delta(t,X)  \overset{\mathbb{P}}{\underset{ T \rightarrow \infty}{\rightarrow}} 0
\end{equation}
when $T \rightarrow \infty$ and when $t$ is distributed according to $\frac{1}{T} \mathbf{1}_{[0,T]}(t) dt$ and $X=(x_{1},\cdots,x_{d})$ is distributed according to $U([- \frac{1}{2}, \frac{1}{2}])^{\otimes^{d}}$. $\overset{\mathbb{P}}{\underset{ T \rightarrow \infty}{\rightarrow}}$ signifies that the convergence occurs in probability.
\end{prop}
\begin{proof}
It is a direct consequence of Proposition $\ref{prop2}$. 
\end{proof}
Because of the independence of the $x_{i}$ between them and with $t$, it is convenient for us to calculate the characteristic function of $\tilde{\Delta}(t,x_{1})$, $\varphi_{\tilde{\Delta}(t,x_{1})}$ (because also the $x_{i}$ are identically distributed). It is the object of the next subsection.
\subsection{Computation of the characteristic function $\varphi_{\tilde{\Delta}(t,x_{1})}$}
The main object of this subsection is to prove the following proposition : 
\begin{prop}
\label{prop6}
For $x$ a real random variable distributed according to the probability measure $\mathbf{1}_{[- \frac{1}{2}, \frac{1}{2}]}(x) dx$ and being independent from $t$, with $t$ being distributed according to the probability measure $\frac{1}{T} \mathbf{1}_{[0,T]}(t) dt$, one has that the characteristic function of $\tilde{\Delta}(t,x)$, $\varphi_{\tilde{\Delta}(t,x)}$ verifies that 
\begin{equation}
\label{eq28}
\varphi_{\tilde{\Delta}(t,x)}(u) = 2\frac{h}{T}  \frac{1-\cos(u)}{u^{2}} + O(\frac{1}{T})
\end{equation}
As a consequence, one has that $\varphi_{\tilde{\Delta}(t,x)}(u) \underset{T \rightarrow \infty}{\rightarrow}  2 \frac{1-\cos(u)}{u^{2}}$. 
\end{prop}
The proof is basically a computation that uses Proposition $\ref{prop4}$ : 
\begin{proof}
According to Proposition $\ref{prop4}$ and because $x$ and $t$ are independent from one another, one has that 
\begin{equation}
\label{eq29}
\varphi_{\tilde{\Delta}(t,x)}(u) = \mathbb{E} \left( \frac{h}{u T}( \sin(u|1-2 \{ x \}|) + \sin(u(1- |1 - 2 \{ x \}|))  + O(\frac{1}{T})\right)
\end{equation}
because of Equation $(\ref{eq24})$ and the $O$ is uniform in $x$ ($h$ can be, and is, chosen so that it does not depend on $x$). \\
Two quick computations give us that : 
\begin{equation}
\label{eq30}
\mathbb{E} \left( \sin(u|1-2 \{ x \}|) \right) = \frac{1-\cos(u)}{u}
\end{equation}
and
\begin{equation}
\label{eq31}
\mathbb{E} \left( \sin(u(1- |1 - 2 \{ x \}|)) \right) = \frac{1-\cos(u)}{u} \textit{.}
\end{equation}
So, with these last two equations and Equation $(\ref{eq29})$, one has that : 
\begin{equation}
\label{eq32}
\varphi_{\tilde{\Delta}(t,x)}(u) = 2\frac{h}{T}  \frac{1-\cos(u)}{u^{2}} + O(\frac{1}{T}) \textit{.}
\end{equation}
By using the fact that $\frac{h}{T} \rightarrow 1$ when $ T \rightarrow \infty$, one gets finally that 
\begin{equation}
\label{eq33}
\varphi_{\tilde{\Delta}(t,x)}(u) \underset{T \rightarrow \infty}{\rightarrow}  2 \frac{1-\cos(u)}{u^{2}} \textit{.}
\end{equation}
\end{proof}
\subsection{Conclusion}
We have now all the necessary tools to prove Theorem $\ref{thm2}$.
\begin{proof}[Proof of Theorem $\ref{thm2}$]
Proposition $\ref{prop5}$ gives us that it is enough to prove $\Delta(t,X)$ converges in law, when $T \rightarrow \infty$. So, we are going to calculate the characteristic function of $\Delta(t,X)$. \\
One has, because $t$ and the $x_{i}$ are independent random variables : 
\begin{equation}
\label{eq34}
\varphi_{\Delta(t,X)}(u) = \prod_{i=1}^{d} \varphi_{\tilde{\Delta}(t,x_{i})}(2^{d-1}u) 
\end{equation}
and because $\Delta(t,X) = 2^{d-1} \sum_{i=1}^{d} \tilde{\Delta}(t,x_{i}) \textit{.}$\\
Furthermore, the $x_{i}$ are identically distributed according to the probability measure $\mathbf{1}_{[-\frac{1}{2}, \frac{1}{2}]}$. So Proposition $\ref{prop6}$ and Equation $(\ref{eq34})$ give us that : 
\begin{equation}
\label{eq35}
\varphi_{\Delta(t,X)}(u) = \prod_{i=1}^{d} 2\frac{h}{T}  \frac{1-\cos(2^{d-1}u)}{(2^{d-1}u)^{2}} + O(\frac{1}{T})  \textit{.}
\end{equation}
So, by making $T$ goes to $\infty$, one has that : 
\begin{equation}
\label{eq36}
\varphi_{\Delta(t,X)}(u) \underset{T \rightarrow \infty}{\rightarrow} (2\frac{1-\cos(2^{d-1}u)}{(2^{d-1}u)^{2}})^{d} \textit{.}
\end{equation}
Theorem $\ref{thm3}$ gives us then the wanted result.
\end{proof}
\bibliographystyle{plain}
\bibliography{bibliographie}
\end{document}